\documentclass[12pt]{amsart}
\usepackage{amsmath} \usepackage{amssymb} \usepackage{anysize}
\marginsize{1.25cm}{1.25cm}{2.5cm}{2.5cm}

\theoremstyle{plain} \newtheorem{theor}{Theorem}[section]
{}
\newtheorem{lemma}{Lemma}[section]

\newtheorem{corol}{Corollary}[section]

\theoremstyle{definition} \newtheorem{defin}{Definition}[section]

\newtheorem{notat}{Notation}[section]
\newtheorem{quest}{Question}[section]

\newtheorem{claim}{Claim}

\theoremstyle{remark}  

\numberwithin{equation}{section}

 \newcommand{\AND}{\text{ and }}
 \newcommand{\OR}{\text{ or }}
 \newcommand{\angbr}[1]{\langle
  #1 \rangle} 
\newcommand{\forces}[2]{\Vdash_{#1} \mbox{``} #2 \mbox{''}}
\newcommand{\notforces}[2]{\not\Vdash_{#1} \mbox{``} #2 \mbox{''}}

\newcommand{\Reals}{{\mathbb R}} \newcommand{\Rationals}{{\mathbb Q}}
\newcommand{\Poset}{{\mathbb P}} \newcommand{\Cantor}{{\mathbb C}}
 
\newcommand{\Naturals}{{\mathbb N}}

\title[Uniformity Invariants]{Comparing the Uniformity Invariants of
  Null Sets for Different Measures} \author[S. Shelah]{Saharon Shelah}
\address{Department of Mathematics, Rutgers University, Hill Center,
  Piscataway, New Jersey, U.S.A. 08854-8019}
\curraddr{Institute of Mathematics\\Hebrew University\\
  Givat Ram, Jerusalem 91904, Israel} \email{shelah@math.rutgers.edu}
\author[J. Stepr\={a}ns]{Juris Stepr\={a}ns} \address{Department of
  Mathematics, York University, 4700 Keele Street, Toronto, Ontario,
  Canada M3J 1P3} \curraddr{} \email{steprans@yorku.ca}
\thanks{The first author's research for this paper was supported by
  the United States-Israel Binational Science Foundation and by grant
  No. NSF-DMS97-04477 of the National Science Foundation of the United
  States. This is number 809 in his personal numbering system.  The
  second author's research for this paper was supported by the Natural
  Sciences and Engineering Research Council of Canada and was
  completed while enjoying the hospitality of Rutger's University. The
  second author expresses his thanks to D. Fremlin and A. Miller for
  enduring a lengthy series of lectures based on an early version of
  this paper, a version even more awkward than the current. }
\begin{document}
\begin{abstract}
  It is shown to be consistent with set theory that the uniformity
  invariant for Lebesgue measure is strictly greater than the
  corresponding invariant for Hausdorff $r$-dimensional measure where
  $0<r<1$.
\end{abstract}
\maketitle
\tableofcontents
\bibliographystyle{plain}
\section{Introduction}
The uniformity invariant for Lebesgue measure is defined to be the
least cardinal of a non-measurable set of reals, or, equivalently, the
least cardinal of a set of reals which is not Lebesgue null. This has
been studied intensively for the past 30 years and much of what is
known can be found in \cite{MR96k:03002} and other standard sources.
Among the well known results about this cardinal invariant of the
continuum is that it can equally well be defined using Lebesgue
measure on $\Reals^n$ without changing the value of the cardinal.
Indeed, equivalent definitions will result by using any Borel
probability measure on any Polish space. However, the question of the
values of uniformity invariants for other, non-$\sigma$-finite Borel
measures is not so easily answered. This paper will deal with the most
familiar class of such measures, the Hausdorff measures for fractional
dimension. Observe that by the previous remarks, the least cardinal of
any non-measurable subset of any $\sigma$-finite set will be the same
as the uniformity invariant for Lebesgue measure.  In other words,
this paper will be concerned with the uniformity invariant of the
ideal of $\sigma$-finite sets with respect to a Hausdorff measure.

It will be shown that given any real number $r$ in the interval
$(0,1)$ it is consistent with set theory that every set of reals of
size $\aleph_1$ is  Lebesgue measurable yet there is a set of reals of
size $\aleph_1$ which is not a null set with respect to
$r$-dimensional Hausdorff measure.  This answers Question~FQ from
D.~Fremlin's list of open questions \cite{fpl}.  However, the
motivation was an attempt to resolve the following question posed by
P.~Komj\'{a}th.
 \begin{quest} \label{Komjath} Suppose that every set of size
   $\aleph_1$ has Lebesgue measure zero.  Does it follow that the
   union of any set of $\aleph_1$ lines in the plane has Lebesgue
   measure zero?
\end{quest}
It is worth noting that this is really a geometric question since
\cite{step.42} provides a negative answer to the version of the
problem in which lines are replaced by their topological and measure
theoretic equivalents.  To see the relationship between
Question~\ref{Komjath} and the topic of this paper consider that it is
easy to find countably many unit squares in the plane such that each
line passes through either the top and bottom or the left and right
sides of at least one of these squares.  It is therefore possible to
focus attention on all lines which pass through the top and bottom of
the unit square. For any such line $L$ there is a pair $(a,b)$ such
that both the points $(a,0) $ and $(1,b)$ belong to $L$. If the
mapping which sends a line $L$ to this pair $(a,b)$ is denoted by
$\beta$ then it is easy to see that $\beta$ is continuous and that if
$S\subseteq [0,1]^2$ is a square of side $\epsilon$ then the union of
$\beta^{-1}S$ has measure $\epsilon$ while $S$ itself has measure
$\epsilon^2$. In other words, the Lebesgue measure of the union of the
lines belonging to $\beta^{-1}X$ is no larger than the 1-dimensional
Hausdorff measure of $X$ for any $X\subseteq [0,1]^2$.  In other
words, a negative answer to Question~\ref{Komjath} would imply that
there is $X\subseteq [0,1]^2$ of size $\aleph_1$ which is not null
with respect to linear Hausdorff measure even though every set of
reals of size $\aleph_1$ is null. The consistency of this will be a
consequence of Corollary~\ref{l:2}.

The proof will rely partially on arguments from \cite{step.36} and
\cite{step.37} in which a single stage of a forcing iteration that
would achieve the desired model was described. The material in
\S\ref{s:cmb} and \S\ref{s:apn} is a reorganized and simplified
version of \S 4, \S 5 and \S 6 of \cite{step.37} which has been
suitably modified for the current context.  The approach taken here
differs from the earlier attempt in that the forcing used is finite
branching rather than infinite branching as in \cite{step.37} and this
allows the use of product forcing along lines similar to those in
\cite{GoShMan93}.  The new ingredient needed is described in
\S\ref{s:ffc}.  The arguments presented in \S\ref{s:tfpo} use ideas
explained in greater detail in \cite{MR1613600}, however familiarity
with that paper is not be required in order to follow the reasoning
presented here.

The logical requirements guiding the organization of the paper may not
have resulted in optimal organization for the purposes of
comprehension. Some readers may prefer to start with \S\ref{s:tfpo}
after having read \S\ref{s:n}, Definition~\ref{d:nxx} and
Definition~\ref{d:1}.
\section{Notation}\label{s:n}
If $X\subseteq \Reals$ and $r \in (0,1)$ then the infimum of all
$\sum_{i=0}^\infty (b_i - a_i)^r$ where $\{(a_i,b_i)\}_{i=0}^\infty$
is a cover of $X$ by intervals of length less than $\epsilon$ is often
denoted by ${\mathcal H}^r_\epsilon(X)$.  The $r$-dimensional
Hausdorff capacity of $X$ is denoted by ${\mathcal H}^r_\infty(X)$ and
is defined to be the infimum of all $\sum_{i=0}^\infty (b_i - a_i)^r$
where $\{(a_i,b_i)\}_{i=0}^\infty$ is a cover of $X$ by arbitrary
intervals. The $r$-dimensional Hausdorff measure of a set $X$ is
denoted by ${\mathcal H}^r(X)$ and, when defined, is equal to its
outer measure $\lim_{\epsilon \to 0}{\mathcal H}_\epsilon^r(X)$.
Since it can be shown\footnote{See Lemma 4.6 in \cite{MR96h:28006} for
  example.} that ${\mathcal H}^r_\infty(X) = 0$ if and only if
${\mathcal H}^r(X) = 0$, in order to establish the main result it
suffices to show that it is consistent with set theory that every set
of reals of size $\aleph_1$ is a Lebesgue null set yet there is a set
of reals of size $\aleph_1$ which is not a null set with respect to
$r$-dimensional Hausdorff capacity.  Actually, it will be more
convenient to work with measures and capacities on the Cantor set and
to replace intervals by dyadic intervals. This amounts to dealing with
net measures as described in \cite{MR96h:28006}. Since $r$-dimensional
net measures have the same null sets as $r$-dimensional Hausdorff
measures this will not be of significance to the results of this
paper.
 \begin{notat}
   For the rest of this paper let $r$ be a fixed real number such that
   $0 < r < 1$.
 \end{notat}
\begin{notat}
  Let $\Cantor$ denote $2^\Naturals$ with the usual product topology.
  Let $\Rationals$ be the tree $\{f~\restriction~n \ : \ f\in \Cantor
  \AND n\in \Naturals\}$.  For $a \in \Rationals$ let $[a] = \{x\in
  \Cantor : a\subseteq x\}$ and for $A\subseteq \Rationals$ let $[A] =
  \bigcup_{a\in A}[a]$.  Let $$\Cantor_\epsilon^{\infty}
  =\left\{A\subseteq \Rationals : \sum_{a\in A}2^{-|a|r} \leq
    \epsilon\right\}$$
  and let $\Cantor_\epsilon =\{A\in
  \Cantor_\epsilon^{\infty} : |A| < \aleph_0\}$.  For consistency of
  notation, let $[\Rationals]^{<\aleph_0}$ be denoted by
  $\Cantor_\infty$. Let $\Cantor_\epsilon^1 = \{A\subseteq \Rationals
  : \ \sum_{t\in A}2^{-|t|} \leq \epsilon\}$.  For $X\subseteq
  \Cantor$ define $\lambda^r_\infty(X) = \inf\left\{\epsilon :
    (\exists A\in \Cantor^{\infty}_\epsilon)\ X\subseteq [A]\right\}
  $.  The usual product measure on $\Cantor$ will be denoted by
  $\lambda$.
\end{notat}

\begin{notat}
  Some notation concerning trees will be established.  By a sequence
  will always be meant a function $f:n \to X$ where $n\in\Naturals$
  and $X$ is some set.  Sequences will occasionally be denoted as
  $n$-tuples $(x_1,x_2,\ldots,x_n)$ and, in particular, singleton
  sequences will be denoted by $(x)$.  If $t$ and $s$ are sequences
  the concatenation of $s$ followed by $t$ will be denoted by $(s, t)$
  . This is consistent with considering $X^k$ to be a set of sequences
  because if $x\in X^k$ and $y \in Y^m$ then $(x,y)\in X^k\times Y^m$.
  If $T$ is a tree then $T$ is a set of sequences closed under
  restriction to initial segments.  If $t\in T$ then ${\mathcal
    S}_T(t) = \{x : (t, (x))\in T\}$. Furthermore, $T\angbr{t} =
  \{s\in T : t\subseteq s \OR s\subseteq t\}$.  If $m\in \Naturals$
  then $T[m] = \{t\in T : |t| = m\}$ and $T[< m] = \{t\in T : |t| <
  m\}$ and $T[\leq m] = \{t\in T : |t| \leq m\}$.
 \end{notat}

\begin{notat}
  For any $X\subseteq \Cantor^m$ and $z\in \Cantor^k$ the set $\{x\in
  \Cantor^{m-k}\ : \ (z, x)\in X\}$ will be denoted by $X(z,\cdot)$.
  Similarly, if $F:\Cantor^m\to X$ is a function and $z\in \Cantor^k$
  then $F(z,\cdot)$ will represent the function
  $F(z,\cdot):\Cantor^{m-k}\to X$ defined by $F(z,\cdot)(x) = F(z,
  x)$.
\end{notat}
\section{Continuous mappings between Lebesgue and Hausdorff measures}\label{s:cmb}

The goal of this section is to exploit the difference between Lebesgue
measure and $r$-dimensional Hausdorff measure. It will be shown that
for any continuous function from the reals to the reals there are
arbitrary small sets in the sense of Lebesgue measure whose pre-image
is as large as desired in the sense of $r$-dimensional Hausdorff
measure. For reasons which will reveal themselves in \S\ref{s:apn} it
will also be necessary to consider similar results for products of the
reals.

\begin{defin}\label{d:sq}
  If $\delta > 0$ is a real number and $X\subseteq \Cantor$ define
  $\square_{\delta}(X)$ by
  $$\square_\delta(X) = \inf\left\{\lambda^r_\infty(X\setminus Z) :
    Z\subseteq \Cantor \AND \lambda(Z) < \delta \right\} .$$
\end{defin}
\begin{lemma} If  \label{jaut}
  $ 1 > \gamma > 0$ and $\epsilon >0$ then there is $\eta > 0$ such
  that for all but finitely many $m \in \Naturals$, any measurable
  $E\subseteq \Cantor$ and any measurable $D\subseteq E$
  $$\square_\eta(D) \geq \square_{\gamma}(E)$$
  provided that $
  \lambda(D\cap [s]) \geq\lambda(E\cap [s])\epsilon$ for each $s \in
  \Rationals[m]$.
\end{lemma}
\begin{proof} 
  Let $\eta >0$ be sufficiently small that
  $$\frac{\gamma }{2 - \gamma} > \frac{4\eta}{\epsilon\gamma}$$
  and
  then let $m \in \Naturals$ be so large that the inequality
  $$
  \frac{ 2^{m(1-r)} \epsilon^r \gamma^r}{4^r} \left(\frac{\gamma
    }{2 - \gamma} - \frac{4\eta} {\epsilon\gamma}\right) >
  \square_{\gamma}(E)$$
  is satisfied.
  
  Suppose that $\lambda(Z) < \eta$ and $A\subseteq \Rationals$ are
  such that $D\subseteq Z\cup [A]$ and $$\sum_{a\in A}2^{-|a|r} <
  \square_\gamma(E) .$$
  Let $A^* = \Rationals[\leq m]\cap A$ and let
  $E^*= E \setminus [A^*]$.  Let $B_0=\{b\in \Rationals[m]\ : \ 
  [b]\cap [A^*] = \emptyset\}$ and note that $E^*\subseteq [B_0]$ and
  that if $b\in B_0$ then $E\cap [b] = E^*\cap [b]$.  Now let $$B_1 =
  \left\{ b \in B_0 \ : \ \lambda(E\cap [b]) >
    \frac{1}{2}\frac{\lambda(E^*)}{2^{m}}\right\}.$$
  Notice that
  $\lambda(E^*) \geq \gamma $ because $\sum_{a\in A^*}2^{-|a|r} <
  \square_\gamma(E)$ and $E\subseteq E^*\cup [A^*]$.  It follows from
  a Fubini argument that
\begin{equation}
  \label{eq:size2}
  |B_1| \geq \frac{2^{m}\lambda(E^*)}{2  -
 \lambda (E^*)} \geq
\frac{2^{m}\gamma}{2 - \gamma} . 
\end{equation}

Next, let $$B_2 = \left\{b\in B_1 \ : \ \lambda(Z\cap [b]) >
  \frac{\epsilon\lambda(E\cap [b])}{2}\right\} $$
and in order to see
that
\begin{equation}
  \label{eq:size1}
  |B_2| \leq \frac{2^{m+2}\eta}{\epsilon\gamma}
\end{equation}
assume the opposite.  Then the following sequence of inequalities
$$\eta > \lambda(Z) \geq \sum_{b\in B_2}\lambda(Z\cap [b]) \geq
\frac{\epsilon}{2}\sum_{b\in B_2}\lambda(E\cap [b]) \geq
\epsilon\sum_{b\in B_2}\frac{\lambda(E^*)}{2^{m+2}} \geq
|B_2|\frac{\epsilon\gamma}{2^{m+2}} \geq \eta$$
yields a
contradiction. (The fourth inequality uses that $B_2 \subseteq B_1$.)

Let $B_3 =B_1\setminus B_2$ and observe that it follows from
Inequalities~\ref{eq:size2} and~\ref{eq:size1} that
\begin{equation}
  \label{eq:size3}
  B_3 \geq \frac{2^{m}\gamma}{2 - \gamma}  -
  \frac{2^{m+2}\eta}{\epsilon\gamma}
=
2^m\left(\frac{\gamma}{2 - \gamma}  - \frac{4\eta}{\epsilon\gamma}\right) .
\end{equation}
It follows that
$$\sum_{a\in A}2^{-|a|r} \geq \sum_{b\in B_3}\left( \sum_{a\in
    A\setminus A^*, a\supseteq b}2^{-|a|r} \right) \geq \sum_{b\in
  B_3}(\lambda(D\cap [b] \setminus Z))^r \geq \sum_{b\in
  B_3}(\lambda(D\cap [b]) -\lambda( Z\cap [b]))^r .
$$
Since $B_3\cap B_2 = \emptyset$ it follows, using the hypothesis on
$D$, that the last term dominates
$$\sum_{b\in B_3}(\lambda(E\cap [b])\epsilon -
(\epsilon/2)\lambda(E\cap [b]))^r \geq
\frac{\epsilon^r}{2^r}\sum_{b\in B_3}(\lambda(E\cap [b]))^r \geq
\frac{\epsilon^r}{2^r}\sum_{b\in
  B_3}\left(\frac{\lambda(E^*)}{2^{m+1}}\right)^r \geq$$
$$
\frac{\epsilon^r}{4^r}|
B_3|\left(\frac{\lambda(E^*)}{2^{m}}\right)^r \geq \frac{ 2^{m(1-r)}
  \epsilon^r \gamma^r}{4^r} \left(\frac{\gamma }{2 - \gamma} -
  \frac{4\eta}{\epsilon\gamma}\right)
$$
and hence $\sum_{a\in A}2^{-|a|r} > \square_\gamma(E)$ which is
impossible.
 \end{proof}
 
 If $X\subseteq \Cantor$ then $F:X\to \Cantor$ will be said to have
 small fibres if and only if $\lambda(F^{-1}\{x\}) = 0$ for each $
 x\in \Cantor$.  The proof of Theorem~\ref{main} and the lemmas
 preceding it will rely on decomposing an arbitrary continuous
 function into a piece that has small fibres and a piece which has
 countable range.

\begin{lemma} Let $0 < \mu <1$ and suppose that $\{X_s\}_{s\in \Rationals}$
  is an indexed family of mutually independent $\{0,1\}$-valued random
  variables, each with mean $\mu$. Suppose \label{manylp} that
  $C\subseteq \Cantor$ is a measurable set and that $F_j: C \to
  \Cantor$ is a measurable function with small fibres for $1\leq j
  \leq n$. For any $\epsilon > 0$ for all but finitely many $m \in
  \Naturals$ the probability that $$\lambda\left (\bigcap_{j=1}^n
    \left( \bigcup_{s\in \Rationals[m],X_s=1}
      F_j^{-1}\left[s\right]\right)\right) > \frac{\mu^n\lambda(C)}{2}
  $$
  is greater than $1 - \epsilon$.
\end{lemma}
\begin{proof}  
  This is Lemma~3.1 in \cite{step.36} or Lemma~6.2 in \cite{step.37}
  except that it is stated here for $\Cantor$ rather than $[0,1]$.
\end{proof}

\begin{lemma}
  Suppose that \label{ions}
\begin{itemize}
\item $E \subseteq \Cantor$ is a measurable set
\item ${\mathcal F}$ is a finite family of measurable functions with
  small fibres from $E$ to $\Cantor$
\item $\gamma > 0$, $\mu > 0$.
\end{itemize}
Then there is $\eta > 0$ such that for any $\epsilon > 0 $ and for any
mutually independent, $\{0,1\}$-valued random variables $\{X_s\}_{s
  \in \Rationals}$ with mean $\mu$ and for all but finitely many $ m
\in \Naturals$ the probability that the inequality
\begin{equation}
  \label{eq:6.3s}
\square_\eta\left (\bigcap_{F\in{\mathcal F}}
F^{-1}\left(
\bigcup_{s\in \Rationals[m],X_s=1}
\left[s\right]\right)\right) 
\geq \square_\gamma(E)  
\end{equation}
holds is greater than $1 - \epsilon$.
 \end{lemma}
\begin{proof} Let $|{\mathcal F}| = n$. Use
  Lemma~\ref{jaut} to choose $\eta > 0$ and an integer $k$ such that
  if $D\subseteq E$ is a measurable set such that for each $s \in
  \Rationals[k]$
  $$
  \lambda(D\cap[s]) \geq \frac{\mu^n}{2}\lambda(E\cap [s])$$
  then
  $\square_\eta(D)\geq \square_{\gamma}(E)$.  Let $\epsilon > 0$. Now
  use Lemma~\ref{manylp} to conclude that for each $t\in\Rationals[k]$
  for all but finitely many $m \in \Naturals$ and any mutually
  independent, $\{0,1\}$-valued random variables $\{X_s\}_{s \in
    \Rationals[m]}$ with mean $\mu$, the probability that
  $$\lambda\left ([t]\cap\bigcap_{F\in{\mathcal F}} F^{-1}\left(
      \bigcup_{s\in \Rationals[m],X_s=1} \left[s\right]\right)\right)
  \geq \frac{\mu^n}{2} \lambda(E\cap [t])
  $$
  is greater than $1 -\epsilon 2^{-k}$. Hence the probability that
  this holds for all $t\in \Rationals[k]$ is greater than $1-\epsilon$
  and so the hypothesis on $k$ guarantees that
  Inequality~\ref{eq:6.3s} holds with at least the same probability.
\end{proof}

\begin{corol}
  Suppose \label{average} that $C$ is a measurable subset of $
  \Cantor^{d+1}$ and ${\mathcal F}$ is a finite family of measurable
  functions from $C$ to $\Cantor$ such that $F(x,\cdot)$ has small
  fibres for each $ x\in \Cantor^d$ and $F\in\mathcal F$. If $\mu > 0$
  and $\gamma > 0$ then there is $\eta > 0$ such that for all
  $\epsilon > 0$ there is some $a\in \Cantor^1_\mu$ such that the
  Lebesgue measure of
  $$\left\{x\in \Cantor^d : \square_{{\eta}}
    \left(\bigcap_{F\in{\mathcal F}}F(x,\cdot)^{-1}a\right) \geq
    \square_\gamma(C(x,\cdot))\right\}$$
  is at least $1 - \epsilon$.
\end{corol}
\begin{proof}  
  This is a standard application of Fubini's Theorem using
  Lemma~\ref{ions} and the Law of Large Numbers.
\end{proof}

\begin{lemma}    Let $E\subseteq \Cantor$ be a measurable set and
  ${\mathcal F}$ a finite family of measurable functions from $E$ to
  $\Cantor$.
 \label{finitelymanyfunctions}
 Then for any $\gamma > 0$ and any $\mu > 0$ there is $\eta >0$ and
 $a\in\Cantor^1_\mu$ such that
 $$\square_\eta \left(\bigcap_{F\in{\mathcal F}}F^{-1}a\right) \geq
 \square_\gamma(E) .$$
\end{lemma} \begin{proof}  
  Let $|{\mathcal F}| = n$.  For each $F\in{\mathcal F}$ let $Y_F =
  \cup\{ F^{-1}\{y\} \ : \ \lambda(F^{-1}\{y\}) > 0\}$ and let
  $\bar{F}$ be defined by
  $$\bar{F}(z) =
\begin{cases}
  F(z) & \text{if } z\in E\setminus Y_F\\
  z  & \text{if } z\in  Y_F .\\
\end{cases}$$
Since each $\bar{F}$ is measurable and has small fibres it is possible
to use Lemma~\ref{ions} to conclude that there is $\eta > 0$ and an
integer $m$ and mutually independent, $\{0,1\}$-valued random
variables $\{X_s\}_{s \in \Rationals[m]}$ with mean $\mu/3$ such that
the probability that
$$\square_{2\eta}\left (\bigcap_{F\in{\mathcal F}} \bar{F}^{-1}\left(
    \bigcup_{s\in \Rationals[m],X_s=1} \left[s\right]\right)\right)
\geq \square_\gamma(E) $$
is greater than $1/2$.

Since the mean of each $X_s$ is $\mu/3$ it is possible to choose $m$
so large that the probability that
$$\lambda\left(\bigcup_{s\in\Rationals[m], X_s=1}
  \left[s\right]\right) < \frac{\mu}{2}$$
is also greater than $1/2$.
Hence there is $a'\in \Cantor^1_{\mu/2}$ such that
$$\square_{2\eta}\left (\bigcap_{F\in{\mathcal F}}
  \bar{F}^{-1}a'\right) \geq \square_\gamma(E) .$$
Now for each
$F\in{\mathcal F}$ choose a finite set $A_F\subseteq \Cantor$ such
that $\lambda(Y_F\setminus F^{-1}A_F) < \eta/n $.  Then let $a\in
\Cantor^1_\mu$ be such that $a' \cup \bigcup_{F\in{\mathcal
    F}}A_F\subseteq a$.  It follows that
$$\square_{\eta}\left(\bigcap_{F\in{\mathcal F}} {F}^{-1}a\right)\geq
\square_{\eta}\left(\bigcap_{F\in{\mathcal F}} \bar{F}^{-1}a \setminus
  \bigcup_{F\in{\mathcal F}}(Y_F\setminus F^{-1}A_F) \right) \geq
\square_{\gamma}(E)$$
because ${F}^{-1}a \supseteq \bar{F}^{-1}a
\setminus (Y_F\setminus F^{-1}A_F)$ for each $F\in{\mathcal F}$.
\end{proof}

For the next definition recall Definition~\ref{d:sq}.
\begin{defin}\label{d:nxx}
  If $\Theta:\Naturals \to \Reals^+$, $\Gamma:\Naturals \to \Reals^+$
  and $X\subseteq \Cantor^d$ then the relation
  $\square_{\Theta,\Gamma}(X)$ will be defined to hold by induction on
  $d$. If $d=1$ then $\square_{\Theta,\Gamma}(X)$ if and only if
  $\square_{\Theta(0)}(X) \geq \Gamma(0)$ whereas, if $d > 1$, then
  $\square_{\Theta,\Gamma}(X)$ holds if and only if
  $\square_{\Theta(d-1)} (\{x \in \Cantor :
  \square_{\Theta,\Gamma}(X(x,\cdot))\}) \geq \Gamma(d-1) $. Define
  $\Theta^+(i) = \Theta(i+1)$ and $\Gamma^+(i) = \Gamma(i+1)$.
\end{defin}
The next lemma establishes that the top-down and bottom-up definitions
of the relation $\square_{\Theta,\Gamma}$ are the same.

\begin{lemma}\label{l:td=bu}
  Let $d\geq 2$.  For any $\Theta:\Naturals \to \Reals^+$,
  $\Gamma:\Naturals \to \Reals^+$ and $X\subseteq \Cantor^{d+1}$ the
  relation $\square_{\Theta,\Gamma}(X)$ holds if and only if
  $\square_{\Theta^+,\Gamma^+}\left(\{z\in \Cantor^d\ : \ 
    \square_{\Theta(0)}(X(z,\cdot)) \geq \Gamma(0)\}\right)$.

\end{lemma}
\begin{proof}
  Proceed by induction on $d$ and observe that the case $d=2$ is
  immediate from Definition~\ref{d:nxx}. Assuming the lemma
  established for $d$ let $X\subseteq \Cantor^{d+2}$. Then the
  following sequence of equivalences establishes the lemma:
  $$\square_{\Theta,\Gamma}(X)$$
  $$\square_{\Theta(d+1)}\left( \left\{ z\in \Cantor^1 \ : \ 
      \square_{\Theta,\Gamma}(X(z,\cdot)) \right\} \right)\geq
  \Gamma(d+1)$$
  $$\square_{\Theta(d+1)}\left( \left\{ z\in \Cantor^1 \ : \ 
      \square_{\Theta^+,\Gamma^+} \left( \{w\in \Cantor^d \ : \ 
        \square_{\Theta(0)}(X(z,\cdot)(w,\cdot))\geq \Gamma(0) \}
      \right) \right\} \right)\geq \Gamma(d+1)$$
  $$\square_{\Theta^+(d)}\left( \left\{ z\in \Cantor^1 \ : \ 
      \square_{\Theta^+,\Gamma^+} \left( \{w\in \Cantor^d \ : \ 
        \square_{\Theta(0)}(X((z,w),\cdot))\geq \Gamma(0) \} \right)
    \right\} \right)\geq \Gamma^+(d)$$
  $$\square_{\Theta^+,\Gamma^+} \left( \left\{ (z, w)\in \Cantor^1
      \times \Cantor^d \ : \ \square_{\Theta(0)}(X((z,w),\cdot))\geq
      \Gamma(0) \right\} \right)
  $$
  $$\square_{\Theta^+,\Gamma^+} \left( \left\{ z\in \Cantor^{d+1} \ :
      \ \square_{\Theta(0)}(X(z,\cdot))\geq \Gamma(0) \right\} \right)
  $$
\end{proof}

\begin{lemma}
  \label{l:tae}
  If $\square_{\Theta,\Gamma}(X) $ holds and $X\subseteq \Cantor^d$
  and $\eta < \Theta(i)$ for each $i< d$ and $A\subseteq \Cantor^d$ is
  such that $\lambda(A) < \eta^d$ then $\square_{\Theta -
    \eta,\Gamma}(X\setminus A) $ holds where $(\Theta - \eta)(i) =
  \Theta(i) - \eta$.
\end{lemma}
\begin{proof}
  Proceed by induction on $d$ using Fubini's Theorem.
\end{proof}

\begin{theor}Let $\Gamma:\Naturals \to
  \Reals^+$ and $\Theta:\Naturals \to \Reals^+$ be functions such that
  $\Theta(i) + \Gamma(i) < 1$ for each $i$. Suppose that $C$ is a
  closed subset of $ \Cantor^{d}$ such that
  $\square_{\Theta,\Gamma}(C)$ holds and that ${\mathcal F}$ is a
  finite family of continuous functions from $C$ to $\Cantor$.  If
  $\mu > 0$ then there is
\label{main} 
$a \in \Cantor^1_\mu$ and $\eta : \Naturals \to \Reals^+$ such that
$\square_{\eta,\Gamma}\left( \bigcap_{F\in {\mathcal F}}F^{-1}a
\right) $.
\end{theor}
\begin{proof}  
  Proceed by induction on $d$, noting that if $d = 1$ then this
  follows from Lemma~\ref{finitelymanyfunctions} by setting $\gamma =
  \Theta(0)$ in that lemma.  So assume that the lemma has been
  established for $d$ and that $C$ is a closed subset of $
  \Cantor^{d+1}$, $\square_{\Theta,\Gamma}(C)$ holds and ${\mathcal
    F}$ is a finite family of continuous functions from $C$ to
  $\Cantor$ and that $\mu > 0$.  For each $F\in\mathcal F$ let $$Y_F =
  \{(x,y)\in \Cantor^d\times\Cantor \ : \ \lambda(\{z\in \Cantor \ : \ 
  F(x,y) = F(x,z)\}) > 0\}$$
  and define $\bar{F}:C \to \Cantor$ by
  $$\bar{F}(x,y) =
\begin{cases}
  F(x,y) & \text{if } (x,y)\notin Y_F\\
  y & \text{otherwise. }
\end{cases}$$
Observe that $\bar{F}(x,\cdot)$ has small fibres for each $x\in
\Cantor^d$.  By Corollary~\ref{average} there is $\delta > 0$ such
that for all $\epsilon > 0$ there is some $a\in \Cantor^1_{\mu/2}$
such that the Lebesgue measure of
$$B(a) = \left\{x\in \Cantor^d : \square_{{\delta}}
  \left(\bigcap_{F\in{\mathcal F}}\bar{F}(x,\cdot)^{-1}a\right) \geq
  \square_{\Theta(0)}(C(x,\cdot))\geq \Gamma(0)\right\}$$
is at least
$1 - \epsilon$.

Since each of the relations $Y_F$ is Borel, it is possible to find a
finite family of functions ${\mathcal G}$ from $\Cantor^d$ to
$\Cantor$ such that
$$\sum_{F\in{\mathcal F}}\int_{x\in \Cantor^d} \lambda\left(
  Y_F(x,\cdot)\setminus F(x,\cdot)^{-1} \{G(x)\}_{G\in{\mathcal G}}
\right)dx < \frac{\delta^{d+1}}{2^{d+1}} .$$
It is then possible to
find a closed set $D\subseteq \Cantor^d$ such that
\begin{itemize}
\item\label{ir:1} $\lambda(D) > \Theta(i+1) + \Gamma(i+1)$ for each $
  i < d$ (using Hypothesis~\ref{eq:fix21})
\item\label{ir:2} $\lambda(D) > 1- \delta^d/2^d$
\item\label{ir:3} $\lambda\left( Y_F(x,\cdot)\setminus F(x,\cdot)^{-1}
    \{G(x)\}_{G\in{\mathcal G}}
  \right) < \delta/2 $ for each $x \in D$
\item\label{ir:4} each $G\in \mathcal G$ is continuous on $D$.
\end{itemize}
Condition~\ref{ir:1} guarantees that $\square_{\Theta^+,\Gamma^+}(D)$
holds.  It is therefore possible to apply the induction hypothesis to
$\Theta^+$, $\Gamma^+$, $\mathcal G$ and $D$ to get $a_0 \in
\Cantor^1_{\mu/2}$ and $\eta:\Naturals \to \Reals^+$ such that
$$\square_{{\eta,\Gamma^+}}\left( D\cap \bigcap_{G\in {\mathcal
      G}}G^{-1}a_0 \right). $$
Let $0 < \epsilon' <\eta(i)$ for $i <
d$.  Now choose $a_1\in \Cantor^1_{\mu/2}$ such that $\lambda(B(a_1))
> 1 - \epsilon'$. Then
$$\square_{{\eta - \epsilon'},\Gamma^+}\left(D\cap B(a_1)\cap
  \bigcap_{G\in {\mathcal G}}G^{-1}a_0 \right)$$
holds by
Lemma~\ref{l:tae} and Condition~\ref{ir:2} and note that by
Condition~\ref{ir:3}
$$\square_{\delta/2} \left(\bigcap_{F\in{\mathcal
      F}}\bar{F}(x,\cdot)^{-1}a_1 \cap Y_F(x,\cdot)\setminus
  F(x,\cdot)^{-1} \{G(x)\}_{G\in{\mathcal G}}
\right) \geq\Gamma(0)$$
for $x\in D\cap B(a_1)$.  Let $\bar{\eta}$ be
defined by
$$\bar{\eta}(i) =
\begin{cases}
  \eta(i-1) - \epsilon' & \text{if }d\geq i >0\\
  \eta(i-1)  & \text{if }d < i \\
  \delta/2 &\text{if } i= 0
\end{cases}$$
and define $Z$ to be the set of all $(x,w) \in \Cantor^{d+1}$ such
that the following three conditions are satisfied:
\begin{equation}
  \label{eq:z1}
x\in D\cap B(a_1)\cap
\bigcap_{G\in {\mathcal G}}G^{-1}a_0  
\end{equation}
\begin{equation}
  \label{eq:z2}
  (\forall F\in{\mathcal F})
\ \bar{F}(x, w)\in a_1
\end{equation}
\begin{equation}
  \label{eq:z3}
  (\forall F\in{\mathcal F})\ F(x, w)\in
Y_F(x,\cdot)\AND F(x,w)\neq G(x)
\end{equation}
It follows from Lemma~\ref{l:td=bu} that
$\square_{\bar{\eta},\Gamma}(Z)$.  Hence it suffices to observe that
if $(x,w)\in Z$ then $F(x,w) \in a_0\cup a_1$. To see this note that
if $F(x,w) = \bar{F}(x,w)$ this is immediate from \ref{eq:z2}.
Otherwise $F(x,w)\in Y_F$ and hence it follows from \ref{eq:z3} that
there is some $G\in\mathcal G$ such that $F(x,w) = G(x)$. From
\ref{eq:z1} it can be concluded that $G(x) \in a_0$.
\end{proof}

\begin{notat}\label{n:1}
  For the rest of the paper, fix a pair of functions $\Gamma:\Naturals
  \to \Reals^+$ and $\Theta:\Naturals \to \Reals^+$ such that
\begin{equation}
  \label{eq:fix1}
\Theta(n) > \Theta(n+1) \AND  \lim_{n\to\infty}\Theta(n) = 0
\end{equation}
\begin{equation}
  \label{eq:fix21}
  (\forall j)\ \Theta(j) + \Gamma(j) < 1
\end{equation}
\begin{equation}
  \label{eq:fix2}
  (\forall j)\ \Theta(j) + \Gamma(j+1) < 1
\end{equation}
\begin{equation}
  \label{eq:fix3}
 \text{the range of }\Gamma\text{ is a dense subset of }(0,1).  
\end{equation}
\end{notat}

\begin{defin}  
  For any $d \in \Naturals$ and for $X\subseteq \Cantor^d$ define
  $\square(X) = \square_{\Theta,\Gamma}(X) $ and define $\square_*(X)$
  to hold if and only if there is some $\eta:\Naturals\to \Reals^+$
  such that $\square_{\eta,\Gamma}(X) $ holds.
\end{defin}

\section{A preliminary norm}
\label{s:apn}
This section will introduce three norms\footnote{This is the
  terminology of \cite{MR1613600}. The exact meaning of the term
  ``norm'' is not required for the arguments used here.} on subsets of
$\Cantor_\infty$ --- these will be denoted $\rho$, $\nu^\infty$ and
$\nu$.  The norms of \cite{MR1613600} typically enjoy some form of
sub-additivity but this will not be the case for any of these three,
at least not explicitly.  Nevertheless, Lemma~\ref{l:5.9s} can be
considered a substitute for this. The only norm used in the definition
of the partial order in \S\ref{s:tfpo} will be $\nu$. The role of the
norm $\rho$ will be to establish a connection between $\nu^\infty$ and
the results of \S\ref{s:cmb}. The norm $\nu^\infty$ is an intermediary
between $\rho $ and $\nu$ and, furthermore, it has the advantage of
allowing the compactness argument of \S\ref{s:ffc} to work.
\begin{notat}
  For any Polish space $X$ let ${\mathcal K}(X)$ denote the spaces of
  compact subsets of $X$ with the Hausdorff metric and let ${\mathcal
    C}(X)$ denote the space of continuous $\Cantor$-valued functions
  with the uniform metric.
\end{notat}
\begin{defin}\label{d:rd}
  Define $\rho$ to be a function from ${\mathcal P}(\Cantor_\infty)$
  to $\Naturals \cup \{\infty\}$ by first defining $\rho^*({\mathcal
    X}) $ to be
  $$\left\{d > 1 \ : \ (\forall C\in {\mathcal K}( \Cantor^d))
    (\forall F\in {\mathcal C}(C)) (\exists x\in {\mathcal X})
    \square(C)\Rightarrow \square_*(F^{-1}x)\right\}$$
  and let
  $\rho({\mathcal X}) = \sup(\rho^*({\mathcal X}))$.  If
  $\rho^*({\mathcal X})=\emptyset$ but ${\mathcal X} \neq \emptyset$
  then define $\rho({\mathcal X}) =0$.
\end{defin}
The following is a direct Corollary of Theorem~\ref{main}.
\begin{corol}\label{c:41}
  If $\mu > 0$ then $\rho(\Cantor^1_\mu) =\infty$.
\end{corol}
\begin{lemma}
  \label{l:5.7.30}
  If $d \geq 1$ then $\{C\subseteq {\mathcal K}(\Cantor^d) \ : \ 
  \square(C)\}$ is a Borel set.
\end{lemma}
\begin{proof}
  It suffices to show $\square(C)$ fails if and only if there is some
  finite subset ${\mathcal C}\subseteq \Rationals^d$ such that
  $$C\subseteq [{\mathcal C}] = \bigcup_{(t_1,t_2,\dots,t_d) \in
    {\mathcal C}} [t_1]\times[t_2]\times \ldots\times[t_d]
  $$
  and $\square([{\mathcal C}])$ also fails.  To prove this proceed
  by induction on $d$. The case $d=1$ is easy. Assuming the result for
  $d$ let $C\subseteq\Cantor^{d+1}$ be compact and suppose that
  $\square(C)$ fails. Using Definition~\ref{d:nxx} there is ${\mathcal
    C}_1\in \Cantor_{\gamma}$ for some $\gamma < \Gamma(d)$ such that
  $\{z\in \Cantor \ : \ \square_{\Theta,\Gamma}(C(x,\cdot))\}\subseteq
  [{\mathcal C}_1]$. For $x\in \Cantor\setminus [{\mathcal C}_1]$ it
  is possible to use the induction hypothesis to find a finite
  ${\mathcal C}_x$ such that $C(x,\cdot) \subseteq [{\mathcal C}_x]$
  and $\square_{\Theta,\Gamma}({\mathcal C}_x)$ fails.  . The
  compactness of $C$ yields an integer $k_x$ such that $C(y,\cdot)
  \subseteq [{\mathcal C}_x]$ for each $y\in \Cantor$ such that
  $y\restriction k_x = x\restriction k_x$.  Let $U_x$ denote this
  neighbourhood and choose a finite $A\subseteq \Cantor$ such that
  $[{\mathcal C}_1]$ and $\{U_a\}_{a\in A}$ cover $\Cantor$.  It is
  then easy to extract the desired finite set from the cover
  $\{[{\mathcal C}_1]\times \Cantor^d\}\cup\{U_a\times[{\mathcal
    C}_a]\}_{a\in A}$.
\end{proof}
\begin{corol}
  \label{c:5.7.30}
  For any $d \geq 1$ the set $\{C\subseteq {\mathcal K}(\Cantor^d) \ :
  \ \square_*(C)\}$ is Borel.
\end{corol}
\begin{proof}
  According to Definition~\ref{d:sq}, $\square_*(C)$ holds if and only
  if there is some $\eta:\Naturals\to \Reals^+$ such that
  $\square_{\eta,\Gamma}(C) $ holds. But notice that $\eta$ can be
  assumed to be a constant function with rational value.
\end{proof}

\begin{notat}
  The notation $\angbr{z, F}$ will be introduced for any $z\in
  \Cantor$ and any function from $F:W \to {\mathcal P}(\Rationals)$ to
  denote $\{w\in W : z \notin [ F(w)]\}$.
\end{notat}

\begin{corol}
  \label{l:5.8.30}
  Suppose that ${\mathcal X}\subseteq \Cantor_\infty$ and $B:{\mathcal
    X} \to {\mathcal P}(\Cantor_\infty)$.  Then $\{z\in \Cantor \ : \ 
  \rho(\angbr{z,B} < j\}$ is analytic for any integer $j$.
\end{corol}
\begin{proof}
  This is immediate from Definition~\ref{d:rd}, Lemma~\ref{l:5.7.30}
  and Corollary~\ref{c:5.7.30}.
\end{proof}

\begin{lemma}
  \label{l:5.9s}
  Suppose that ${\mathcal X}\subseteq \Cantor_\infty$ and
  $\rho({\mathcal X}) = j\geq 2$.  Then
  $$\square_{\Theta(j-2)}\left(\left\{ z\in \Cantor \ : \ 
      \rho(\angbr{z,F}) < j - 1\right\}\right) < \Gamma(j-1) $$
  for
  each function $F$ from ${\mathcal X}$ to
  $\Cantor^\infty_{\Gamma(j-1)/2}$.
\end{lemma}
\begin{proof}
  Let $S = \left\{ z\in \Cantor \ : \ \rho(\angbr{z,F}) < j -
    1\right\}$ and assume that the lemma fails. In other words,
  $\square_{\Theta(j-2)}(S) \geq \Gamma(j-1)$.  Let
  $$E\subseteq \Cantor \times {\mathcal K}(\Cantor^{j-1})\times
  {\mathcal C}(\Cantor^{j-1})$$
  be the set of all triples $(z, C, f)$
  such that
\begin{equation}
  \label{eq:m301}
  \square(C)
\end{equation}
\begin{equation}
  \label{eq:m302}
  (\forall x \in \angbr{z,F})\ \neg\square_*(C\cap f^{-1}x)
\end{equation}
It follows from Definition~\ref{d:rd} that $S$ is contained in the
domain of $E$.  From Lemma~\ref{l:5.7.30} and Corollary~\ref{c:5.7.30}
it follows that $E$ is a Borel set.  From Corollary~\ref{l:5.8.30} it
follows that $S$ is measurable and so it is possible to use the
von~Neumann Selection Theorem and Egeroff's Theorem to find a closed
set $S^*$ in the domain of $E$ and a continuous function $T\subseteq
E$ with domain $S^*$ such that $\lambda(S\setminus S^*) < \Theta(j-2)
- \Theta(j-1)$. Let $T(s) = (C_s, f_s)$ and observe that the
continuity of $T$ guarantees that $C^* = \bigcup_{s\in S^*}\{s\}\times
C_s$ is a closed set. Moreover, $\square_{\Theta(j-1)}(S^*) \geq
\square_{\Theta(j-2)}(S) \geq \Gamma(j-1)$. A calculation using
Definition~\ref{d:nxx} reveals that $\square(C^*) $ holds.  Let $g$ be
defined on $C^*$ by $g(c_1, c_2, \ldots, c_j) =
f_{c_1}(c_2,c_3,\ldots,c_j)$ recalling that $j\geq 2$.

Since $\rho({\mathcal X}) \geq j$ it is possible to find $x\in
{\mathcal X}$ such that $\square_*(g^{-1}x)$.  Let $$W= \{z\in S^* \ :
\ \square_*(g^{-1}x(z,\cdot))\} =\{z\in S^* \ : \ 
\square_*(g(z,\cdot)^{-1}x)\}
$$
and note that it follows that $ \lambda^r_\infty(W) \geq
\Gamma(j-1)$.  Since $F(x)\in \Cantor_{\Gamma(j-1)/2}$ it is possible
to choose $w\in W\setminus F(x)$. Then $x\in \angbr{w,F}$ and so, by
Condition~\ref{eq:m302} in the definition of $E$, it follows that
$\square_*( g(w,\cdot)^{-1}x)$ fails since $g(w,\cdot) = f_w$. This
contradicts that $w\in W$.
\end{proof}
\begin{defin}\label{d:1}
  The norms $\nu$ and $\nu^\infty $ will be defined for the subsets of
  $\Cantor_\infty$ by first using induction to define an associated
  sequence of sets:
  $${\mathcal N}_0 = {\mathcal N}_0^\infty = \{{\mathcal
    X}\subseteq\Cantor_\infty \ : \ {\mathcal X}\neq \emptyset\}$$
  $${\mathcal N}_1 = {\mathcal N}_1^\infty =\{{\mathcal
    X}\subseteq\Cantor_\infty \ : \ [\cup {\mathcal X}]=\Cantor\}$$
\begin{equation}
  \label{eq:mar25.2}
{\mathcal N}_{j+1} = \{{\mathcal X}\subseteq\Cantor_\infty \ : \ 
(\forall F:{\mathcal X} \to \Cantor_{\Gamma(j)/2}) \
\square_{\Theta(j-1)}\left(\{z \in \Cantor :
\angbr{z,F}\notin {\mathcal N}_j)\}\right) < \Gamma(j)
\}
\end{equation}
\begin{equation}
  \label{eq:mar25.1}
{\mathcal N}_{j+1}^\infty =\{{\mathcal X}\subseteq\Cantor_\infty \ : \ 
(\forall F:{\mathcal X} \to \Cantor_{\Gamma(j)/2}) \
\square_{\Theta(j-1)}\left(\{z \in \Cantor :
\angbr{z,F}\notin {\mathcal N}_j^\infty)\}\right) < \Gamma(j)
\} .\end{equation}
Then define
$\nu({\mathcal X})$ to be the supremum of all $j$ such that 
${\mathcal X}\in {\mathcal N}_j$ and
$\nu^\infty({\mathcal X})$ to be the supremum of all $j$ such that 
${\mathcal X}\in {\mathcal N}^\infty_j$.
\end{defin}

It must be noted that ${\mathcal N}_{j}\supseteq {\mathcal N}_{j+1}$
and ${\mathcal N}^\infty_{j}\supseteq {\mathcal N}^\infty_{j+1}$ for
each integer $j$ and hence the supremum in Definition~\ref{d:1} is
taken over an initial segment of the integers. This will be used
implicitly in what follows.

\begin{corol}
  \label{c:5.2}
  If ${\mathcal X}\subseteq \Cantor_\infty$ then $\rho({\mathcal
    X})\leq \nu^\infty ({\mathcal X})$.
\end{corol}
\begin{proof}
  Proceed by induction on $j = \rho({\mathcal X})$. The case $j=0$ is
  trivial but the case $j=1$ is less so. To see that $\rho({\mathcal
    X}) = 1$ implies that ${\mathcal X}$ covers $\Cantor$ let $z\in
  \Cantor$. Note that $\square(\Cantor)$ holds by
  Condition~\ref{eq:fix2} of Notation~\ref{n:1}.  Letting $F$ be the
  function on $\Cantor$ with constant value $z$ it follows from the
  definition of $\rho({\mathcal X}) = 1$ that there is some $x\in
  {\mathcal X}$ such that $\square_*(F^{-1}x) $ holds.  In particular,
  $F^{-1}x \neq \emptyset$ and so if $w\in F^{-1}x$ then $z=F(w)\in
  x$.
  
  Therefore it can be assumed that $2 \leq j+1 = \rho({\mathcal X})$
  and that the lemma holds for $j$. In order to show that
  $\nu^\infty({\mathcal X}) \geq j+1$ let $F:{\mathcal X} \to
  \Cantor^\infty_{\Gamma(j)/2}$.  By Lemma~\ref{l:5.9s} it must be
  that
  $$\square_{\Theta(j-1)}\left(\left\{ z\in \Cantor \ : \ 
      \rho(\angbr{z,F}) < j\right\}\right) < \Gamma(j) . $$
  From the
  induction hypothesis it follows that if $\nu^\infty(\angbr{z,F}) <
  j$ then $\rho(\angbr{z,F}) < j$.  Hence
  $$\square_{\Theta(j-1)}\left(\left\{ z\in \Cantor \ : \ 
      \nu^\infty(\angbr{z,F}) < j\right\}\right) < \Gamma(j)$$
  and
  this establishes that $\nu^\infty({\mathcal X}) \geq j+1$.
\end{proof}
\section{Finding finite sets with large norm}\label{s:ffc}
It will be shown that for any $j$ there is a finite ${\mathcal X}$
such that $\nu^\infty({\mathcal X}) \geq j $. This will establish that
the ${\mathfrak C}(k)$ required in the definition of the partial order
$\Poset$ in \S\ref{s:tfpo} actually do exist. Each of the next lemmas
is a step towards this goal.
\begin{lemma}\label{l:mnutau}
  If ${\mathcal X}\subseteq\Cantor_\infty$ then $\nu^\infty({\mathcal
    X}) \leq \nu({\mathcal X}) $.
\end{lemma}
\begin{proof}
  This follows from the definitions by an argument using induction on
  $j$ to show that if $\nu^\infty({\mathcal X}) \geq j+1 $ then
  $\nu({\mathcal X}) \geq j+1$.
\end{proof}
\begin{lemma}\label{l:nuvstau}
  If $\{{\mathcal A}_n\}_{n=0}^\infty$ is an increasing sequence of
  finite subsets of $\Cantor_\infty$ then
  $$\nu^\infty\left(\bigcup_{n=0}^\infty {\mathcal A}_n\right) \ = \ 
  \lim_{n\to\infty}\nu^\infty({\mathcal A}_n) .$$
\end{lemma}
\begin{proof}
  Proceed by induction on $j$ to show that if $\nu^\infty({\mathcal
    A}_n) < j $ for each $n$ then
  $\nu^\infty\left(\bigcup_{n=0}^\infty {\mathcal A}_n\right) < j$.
  For $j = 0 $ this trivial and if $j=1$ this is simply a restatement
  of the compactness of $\Cantor$. Therefore assume that $j\geq 1$,
  that the lemma is true for $j$, that $\nu^\infty({\mathcal A}_n)
  \leq j $ for each $n$ yet $\nu^\infty\left(\bigcup_{n=0}^\infty
    {\mathcal A}_n\right) > j$.  Let $F_n :{\mathcal A}_n \to
  \Cantor^{\infty}_{\Gamma(j)/2}$ witness that $\nu^\infty({\mathcal
    A}_n)\not\geq j+1 $.  In other words, using
  Equality~\ref{eq:mar25.1} of Definition~\ref{d:1},
  $\square_{\Theta(j-1)}(S_n) \geq \Gamma(j)$ where $S_n = \{z\in
  \Cantor : \nu^\infty(\angbr{z,F_n}) < j\}$.
\begin{claim}
  If $\{\delta_n\}_{n=0}^\infty$ is a sequence of positive reals and
  and $\{A_n\}_{n=0}^\infty$ is a sequence of elements of
  $\Cantor^{\infty}_{1}$ then there are two increasing sequences of
  integers $\{k_n\}_{n=0}^\infty$ and $\{m_n\}_{n=0}^\infty$ such
  that, letting $$D_n=A_{m_n}\cap \Rationals[<k_n]$$
  the following
  hold:
\begin{equation}
  \label{eq:281}
  \text{ if }i\geq n \text{ then } 
D_n=D_i\cap \Rationals[<{k_n}]
\end{equation}
\begin{equation}
  \label{eq:282}
  \lambda\left(\bigcup_{i=n}^\infty[A_{m_i}\setminus D_n]\right) < \delta_n.
\end{equation}
Moreover, the increasing sequence $\{k_n\}_{n=0}^\infty$ can be chosen
from any given infinite set $K$.
\end{claim}
\begin{proof}
  Let $K$ be given and, using the fact that $0 < r < 1$, let
  $\{k_n\}_{n=0}^\infty\subseteq K$ be such that for each $n$
  $$\sum_{i=n}^\infty2^{-k_i(1-r)} \leq \delta_n/2 $$
  and then choose
  the sequence $\{m_n\}_{n=0}^\infty$ such that
  Conclusion~\ref{eq:281} holds.  Then Conclusion~\ref{eq:282} follows
  from the following inequalities
  
  $$
  \lambda\left(\bigcup_{i=n}^\infty[A_{m_i}\setminus D_n]\right)
  \leq \lambda\left(\bigcup_{i=n}^\infty[A_{m_i}\setminus D_i]\cup
    [D_{i+1}\setminus D_i]\right) \leq \sum_{i=n}^\infty
  \lambda([A_{m_i}\setminus D_i]) + \lambda([D_{i+1}\setminus D_i])$$
  $$\leq \sum_{i=n}^\infty\left( \sum_{t\in A_{m_i}\setminus
      D_i}2^{-|t|} + \sum_{t\in D_{i+1}\setminus D_i}2^{-|t|}\right)
  \leq \sum_{i=n}^\infty\left( \sum_{t\in A_{m_i}\setminus
      D_i}2^{-|t|r}2^{-|t|(1-r)} + \sum_{t\in D_{i+1}\setminus
      D_i}2^{-|t|r}2^{-|t|(1-r)} \right)$$
  $$\leq \sum_{i=n}^\infty\left( 2^{-k_i(1-r)} \left(\sum_{t\in
        A_{m_i}\setminus D_i}2^{-|t|r} + \sum_{t\in D_{i+1}\setminus
        D_i}2^{-|t|r} \right)\right) \leq \sum_{i=n}^\infty
  2^{-k_i(1-r)}\left(1 + 1\right) \leq \delta_n .$$
\end{proof}

Before continuing, let $\{a_i\}_{i=0}^\infty$ enumerate
$\bigcup_{n=0}^\infty {\mathcal A}_n$ and, without loss of generality,
assume that ${\mathcal A}_n = \{a_\ell\}_{\ell = 0}^n$.  Using the
claim and its final clause and a diagonalization argument, find two
increasing sequences of integers $\{k_n\}_{n=0}^\infty$ and
$\{m_n\}_{n=0}^\infty$ such that for each $i$ and each $n\geq i$,
letting $F_{n,i}= F_{m_n}(a_i)\cap \Rationals[<k_n]$,
Conclusion~\ref{eq:281} and Conclusion~\ref{eq:282} of the claim hold
for a summable sequence of $\delta$.  To be more precise,
$F_{{n+1},i}$ is an end extension of $F_{{n},i}$ and
$$\lambda\left(\bigcup_{j=n}^\infty\left[F_{m_j}(a_i)\setminus
    F_{n,i}\right]\right) < 2^{-i-n} $$
for any $n\geq i$.  Define
$F(a_i) = \bigcup_{n=i}^\infty F_{n,i}$ and note that $F(a_i) \in
\Cantor^{\infty}_{\Gamma(j)/2}$ since $F_{n,i} \in
\Cantor_{\Gamma(j)/2}$ for each $n$ and the $F_{n,i}$ are increasing
with respect to $n$.

Now let $S= \{z\in \Cantor : \nu^\infty(\angbr{z,F}) < j\}$. Because
$$\nu^\infty\left(\bigcup_{n=0}^\infty {\mathcal A}_n\right) \geq
j+1$$
it follows that $\square_{\Theta(j-1)}(S) < \Gamma(j)$ and,
hence, it is possible to choose $Y$ such that
$\lambda^r_\infty(S\setminus Y) < \Gamma(j)$ and $\lambda(Y) <
\Theta(j-1)$.  Choose $M$ so large that letting $$W =
\bigcup_{i=0}^\infty\left(\bigcup_{n=\max(M,i)}^\infty \left[
    F_{m_n}(a_i)\setminus F_{{\max(M,i)},i}\right]\right) $$
it
follows that $\lambda(W) < \Theta(j-1) - \lambda(Y)$.

Now define $F^*(a_i)$ by
$$F^*(a_i) =
\begin{cases}
  F_{M,i} & \text{if } i\leq M\\
  F_{i,i} & \text{if } i > M
\end{cases}$$
and note that $F^*(a_i)\subseteq F(a_i)$.  Now let $S^*_m =
\{z\in\Cantor : \nu^\infty(\angbr{z,F^*\restriction {\mathcal A}_m}) <
j\}$ and note that $S^*_m\supseteq S^*_{m+1}$ and they are all
compact.  Indeed,
$$S^*_m = \bigcup_{{\mathcal B}\subseteq {\mathcal A}_m,
  \nu^\infty({\mathcal B}) < j} \left( \bigcap_{ a\in {\mathcal
      A}_m\setminus {\mathcal B}} [F^*(a)]\right)$$
and each
$[F^*(a)]$ is compact.

\begin{claim}
  If $z\notin W$ and $i\geq M$ then $\angbr{z,F^*\restriction{\mathcal
      A}_{m_i}}\subseteq \angbr{z,F_{m_i}}$.
\end{claim}
\begin{proof}
  Let $a_\ell \notin \angbr{z,F_{m_i}}$ and assume that $\ell \leq
  m_i$ since otherwise it is immediate that
  $a_\ell\notin\angbr{z,F^*\restriction {\mathcal A}_{m_i}}$.  There
  are two cases to consider. First assume that $\ell \leq M$.  Then
  $z\in [F_{m_i}(a_\ell)]$ but since $z\notin W$ and $i\geq M$ it
  follows that $z \notin [ F_{m_i}(a_\ell)\setminus F_{{M},\ell}]$.
  Hence $z \in [F_{M,\ell}] = [F^*(a_\ell)]$. In other words,
  $a_\ell\notin\angbr{z,F^*\restriction {\mathcal A}_{m_i}}$.  If
  $\ell > M$ a similar argument works.
\end{proof}

It follows from the claim that if $z\notin W$ and $i \geq M$ then
$\nu^\infty(\angbr{z,F^*\restriction{\mathcal A}_{m_i}}) \leq
\nu^\infty(\angbr{z,F_{m_i}})$. Hence, if $i \geq M$ and $z\in
S_{m_i}\setminus W$ then $\nu^\infty(\angbr{z,F^*\restriction{\mathcal
    A}_{m_i}}) < j$. In other words, $S_{m_i}\setminus W\subseteq
S_{m_i}^*\setminus W$.  Therefore, $\square_{\Theta(j-1) -
  \lambda(W)}(S^*_{m_i}\setminus W)\geq \square_{\Theta(j-1) -
  \lambda(W)}(S_{m_i}\setminus W) \geq \Gamma(j)$.  Since each of the
$S_{m_i}^*\setminus W$ are compact it follows that
$$\square_{\Theta(j-1) - \lambda(W)}\left(\bigcap_{i=M}^\infty
  S_{m_i}^*\setminus W\right) \geq \Gamma(j) . $$
Let $S^* =
\bigcap_{i=M}^\infty S^*_i$.  Since $\lambda(Y) < \Theta(j-1) -
\lambda(W)$ it follows that $\lambda^r_\infty\left( S^*\setminus
  (W\cup Y)\right) \geq \Gamma(j)$.

Since $\lambda^r_\infty(S\setminus Y) < \Gamma(j)$ it is possible to
select $z\in S^*\setminus (S\cup W \cup Y)$.  Then
$\nu^\infty(\angbr{z,F^*\restriction{\mathcal A}_{m_i}}) < j$ for each
$i\geq M$ and so the induction hypothesis guarantees that
$$\nu^\infty\left(\bigcup_{i=M}^\infty\angbr{z,F^*\restriction{\mathcal
      A}_{m_i}}\right) < j .$$
But
$$\bigcup_{i=M}^\infty\angbr{ z,F^*\restriction{\mathcal A}_{m_i}} =
\angbr{z,\bigcup_{i=M}^\infty F^*\restriction{\mathcal A}_{m_i}} =
\angbr{z,F^*}$$
and so $\nu^\infty(\angbr{z,F^*}) < j$. Since $F^*(a)
\subseteq F(a)$ for each $a$ it is immediate that
$\angbr{z,F^*}\supseteq \angbr{z,F}$ and so $\nu^\infty(\angbr{z,F}) <
j$.  This contradicts that $z\notin S$.
\end{proof}
\begin{corol}\label{l:sb}
  For any $j \in \Naturals$ and $\mu > 0 $ there is a finite set
  ${\mathcal X}\subseteq \Cantor^1_\mu$ such that $\nu({\mathcal X})
  \geq j$.
\end{corol}
\begin{proof}
  Combine Corollary~\ref{c:41} and Corollary~\ref{c:5.2} to conclude
  that $\nu^\infty(\Cantor^1_\mu)=\infty$. Then use
  Lemma~\ref{l:nuvstau} to find a finite ${\mathcal X}\subseteq
  \Cantor^1_\mu$ such that $\nu^\infty({\mathcal X})\geq j$. Finally
  apply Lemma~\ref{l:mnutau}.
\end{proof}
\section{The forcing partial order}\label{s:tfpo}
Using Corollary~\ref{l:sb} let ${\mathfrak C}(n)$ be a finite subset
of $\Cantor^1_{2^{-n}}$ such that $\nu({\mathfrak C}(n)) \geq n$ for
each $n\in \Naturals$.  Recalling the notation concerning trees in
\S\ref{s:n}, let $\Poset$ consist of all trees $T$ such that:
\begin{equation}\label{eq:1}
(\forall t \in T)(\forall i < |t|)\ t(i) \in {\mathfrak C}(i)
\end{equation}
\begin{equation}\label{eq:2}
(\forall k\in \omega)(\forall t\in T)(\exists s\in T)\ t\subseteq s
\AND
\nu({\mathcal S}_T(s)) > k
\end{equation}
and let this be ordered under inclusion.  The methods of
\cite{MR1613600} can be used to establish that $\Poset$ is an
$\omega^\omega$-bounding proper partial order.
\begin{lemma}\label{l:1}
  Let $V\subseteq W$ be models of set theory and suppose that
  $G\subseteq \Poset\cap V$ is generic over $W$. Then $W[G]\models
  \lambda(W\cap [0,1]) = 0$.
\end{lemma}
\begin{proof}
  If $G\subseteq \Poset\cap V$ is generic over $W$ then let $B_G\in
  \prod_{n=0}^\infty {\mathfrak C}(n)$ be the generic branch
  determined by $G$. Note that $\lambda([B_G(j)]) < 2^{-j}$ for each
  $j$ and so $\lambda(\bigcup_{j=0}^\infty [B_G(j)]) < \infty$. Also,
  for every $T\in \Poset$ there is some $t\in T$ such that $V\models
  \nu({\mathcal S}_T(t)) \geq 1$. Since $\nu({\mathcal S}_T(t)) \geq
  1$ is equivalent to $\cup{\mathcal S}_T(t) = \Cantor$, and this is
  absolute, it follows from genericity that $\Cantor\cap W\subseteq
  \bigcup_{j=m}^\infty [B_G(j)]$ for every $m$.
\end{proof}
\begin{defin}
  Let $\Poset^\kappa$ be the countable support product of $\kappa$
  copies of $\Poset$.
\end{defin}
\begin{corol}\label{c:main}
  If $\kappa \geq \aleph_1$ and $G\subseteq \Poset^\kappa$ is a filter
  generic over $V$ then in $V[G]$ every set of reals of size less than
  $\kappa$ has Lebesgue measure 0.
\end{corol}
In light of Corollary~\ref{c:main} the goal now is to establish that
if $\Poset^\kappa$ is the countable support product of $\kappa$ copies
of $\Poset$ and $G\subseteq \Poset^\kappa$ is a filter generic over
$V$ then $\lambda^r_\infty(\Cantor\cap V) = 1$ in $V[G]$.

\begin{lemma}\label{l:0}
  Suppose that $p\in \Poset^\kappa $ and
  $$p\forces{\Poset^\kappa}{C\in \Cantor^{\infty}_\delta } $$
  and
  $\delta < 1$.  Then there are $A\subseteq \Cantor$ and $B\subseteq
  \Cantor$ such that $\lambda^r_\infty(A) + \lambda(B) < 1$ and $q\leq
  p$ such that $q\notforces{\Poset^\kappa}{z\in [C]}$ for each $z \in
  \Cantor\setminus (A \cup B)$.
\end{lemma}
\begin{proof}
  Let $1 > \delta_1 > \delta$ and choose a monotonically increasing
  function $\delta_2:\Naturals \to \Reals$ such that $\delta_1 +
  \lim_{n\to\infty}\delta_2(n) < 1$. For later use, define $k(i) =
  \prod_{j=0}^i|{\mathfrak C}(j)|$.  The proof is based on a standard
  fusion argument. Induction will be used to construct $p_i$, $S_i$,
  $\epsilon_i$, $a_i$, $b_i$ and $L_i$ satisfying the following
  conditions:
\begin{enumerate}
\item $p_i\in \Poset^\kappa$ and $p_0 \leq p$
\item $p_{i+1}\leq p_i$
\item $S_i :\kappa \to \Naturals$ is 0 at all but finitely many
  ordinals --- these will be denoted by $D(S_i)$
\item $S_i(\sigma) = S_{i+1}(\sigma)$ for all but one $\sigma
  \in\kappa$
\item $S_{i+1}(\sigma) \leq S_i(\sigma) + 1$ for all $\sigma
  \in\kappa$
\item $L_i\in \Naturals$ and $L_i < L_{i+1}$
\item \label{co:7-1}$p_{i+1}(\sigma)[L_i] = p_i(\sigma)[L_i]$ for each
  $\sigma\in D(S_i)$
\item \label{co:7} if $S_i(\sigma) > 0$ there is a maximal antichain
  ${\mathcal C}_{i,\sigma}\subseteq p_i(\sigma)[<L_i]$ such that
  $\nu({\mathcal S}_{p_i(\sigma_j)}(s)) \geq S_i(\sigma)$ for each
  $s\in {\mathcal C}_{i,\sigma}$
\item \label{co:99}for each $\tau\in \prod_{\sigma \in
    D(S_i)}p_i(\sigma)[L_i]$ there is some finite $C_{i,\tau}$ and
  $k_{i,\tau}\in \Naturals$ such that
  $$p_i\angbr{\tau}\forces{\Poset^\kappa}{ \check{C}_{i,\tau}= C\cap
    \Rationals[<k_{i,\tau}] \AND C\setminus \check{C}_{i,\tau}
    \in\Cantor^{\infty}_{\epsilon_i}}$$
  where $p_i\angbr{\tau}$ is
  defined by $p_i\angbr{\tau}(\alpha) =
  p_i(\alpha)\angbr{\tau(\alpha)}$
\item \label{co:10} $a_i\in \Cantor_{\alpha_i}$ for some $\alpha_i <
  \delta_1$ and $a_i\subseteq a_{i+1}$
\item \label{co:101} $\lambda([b_i]) < \delta_2(i)$ and $b_i\subseteq
  b_{i+1}$
\item \label{co:1300}for each $ z\in \Cantor\setminus [a_i\cup b_i]$
  there is $T_{\sigma,z,i}\in \Poset$ for each $\sigma\in D(S_i)$ such
  that
 \begin{enumerate}
 \item $T_{\sigma,z,i}\subseteq p_i({\sigma})[\leq L_i]$
 \item \label{co:11b}if $s\in T_{\sigma,z,i}[<L_i]\cap {\mathcal
     C}_{j,\sigma}$ then $\nu({\mathcal S}_{T_{\sigma,z,i}}(s)) \geq
   S_j(\sigma)$
 \item \label{co:11c}$z\notin [C_{i,\tau}]$ for any $\tau\in
   \prod_{\sigma \in D(S_i)}T_{\sigma,z,i}[L_i]$
 \item \label{co:11d}$T_{\sigma,z,i}[L_i] = T_{\sigma,z,i+1}[L_i]$
 \end{enumerate}
\item\label{co:sm} $2\epsilon_i k(L_i)^{2i} < \delta_1 - \alpha_i$
\item \label{co:17} $\lim_{n\to\infty}S_n(\sigma) = \infty$ for each
  element $\sigma$ of the domain of some $p_i$
\end{enumerate}
Assuming this can be done, let $A = \bigcup_{i=0}^\infty [a_i]$ and $B
= \bigcup_{i=0}^\infty [b_i]$ and define
$$q(\sigma) = \bigcup_{i=0}^\infty p_i(\sigma)[L_i]$$
for each $\sigma
\in\kappa$ which is in the domain of some $p_i$.  It follows from
Conditions~\ref{co:7} and~\ref{co:7-1} that $q\in \Poset^\kappa$.
Then if $z\in \Cantor\setminus [A\cup B]$ let
$$q_z(\sigma) = \bigcup_{j=0 }^\infty T_{\sigma,z,j}$$
and note that
it follows from Condition~\ref{co:11b} that $q_z\in \Poset^\kappa$.
\relax From Conditions~\ref{co:99} and~\ref{co:11c} it follows that
$q_z\forces{\Poset^\kappa}{z\notin C}$.

To see that the induction can be carried out, let $i$ be given and
suppose that $p_i$, $S_i$, $\epsilon_i$, $a_i$, $b_i$ and $L_i$
satisfying the induction hypothesis have been chosen. Choose
$\bar{\sigma}$ according to some scheme which will satisfy
Condition~\ref{co:17} and define
$$S_{i+1}(\sigma) =
\begin{cases}
  S_i(\sigma) &\text{if } \sigma \neq \bar{\sigma}\\
  S_i(\sigma)+1 &\text{if } \sigma = \bar{\sigma}.
\end{cases}$$
Let $\bar{p}_0\leq p_i$ be such that $\bar{p}_0 = p_i$ if
$S_i(\bar{\sigma})>0$ and, otherwise, $\bar{p}_0(\sigma) =
p_i(\sigma)$ for $\sigma \neq \bar{\sigma}$ and
$|\bar{p}_0(\bar{\sigma})[L_i] | = 1$.  Let $m \geq
S_{i+1}(\bar{\sigma})$ be such that
\begin{equation}
  \label{eq:r1}
 \frac{\Gamma(m)}{2} >  \epsilon_i k(L_i)^{i}
\end{equation}
\begin{equation}
  \label{eq:r2}
  \Gamma(m)k(L_i) < \delta_1 - \alpha_i
\end{equation}
\begin{equation}
  \label{eq:r3}
  \Theta(m-1) < \frac{\delta_2(i+1) - \delta_2(i)}{k(L_i)}
\end{equation}
and for each $t\in \bar{p}_0(\bar{\sigma})[L_i]$ there is some $t^*\in
\bar{p}_{0}(\bar{\sigma})$ such that $t\subseteq t^*$ and such that
\begin{equation}
  \label{eq:r4}
 \nu({\mathcal S}_{\bar{p}_0(\bar{\sigma})}(t^*)) \geq m + 1 . 
\end{equation}
Let $L_{i+1}> L_i$ be so large that $|t^*| < L_{i+1}$ for each $t\in
\bar{p}_0(\bar{\sigma})[L_i]$.  Then let $\bar{p}_1\leq \bar{p}_0$ be
such that
 \begin{itemize}
 \item ${\mathcal C}_{i+1,\bar{\sigma}} = \{t^* \ : \ 
   t\in\bar{p}_0(\bar{\sigma} )[L_i]\}$ is a maximal antichain in
   $\bar{p}_1(\bar{\sigma})$
 \item if $t \in \bar{p}_1(\bar{\sigma})$ and $L_i \leq |t|\leq
   L_{i+1}$ then either
  \begin{itemize}
  \item $|{\mathcal S}_{\bar{p}_1(\bar{\sigma})}(t)| = 1$ or
  \item ${\mathcal S}_{\bar{p}_1(\bar{\sigma})}(t) = {\mathcal
      S}_{\bar{p}_0(\bar{\sigma})}(t) $
  \end{itemize}
\item $\bar{p}_1(\sigma) = \bar{p}_0(\sigma)$ for $\sigma \neq
  \bar{\sigma}$.
 \end{itemize}
 Let ${\mathcal C}_{i+1,{\sigma}} = {\mathcal C}_{i,{\sigma}} $ if
 $\sigma \neq \bar{\sigma}$.  Let $\bar{p}_2\leq \bar{p}_1$ be such
 that $\bar{p}_2(\sigma) = \bar{p}_1(\sigma)$ for $\sigma \notin
 D(S_i)\setminus\{\bar{\sigma}\}$ and if $\sigma \in
 D(S_i)\setminus\{\bar{\sigma}\}$ and $t\in \bar{p}_2(\sigma)$ and
 $L_i \leq |t| \leq L_{i+1}$ then $|{\mathcal
   S}_{\bar{p}_2(\sigma)}(t)| = 1$.  Let $\epsilon_{i+1}>0 $ be so
 small that $2\epsilon_{i+1}k(L_{i+1})^{2(i+1)} < \delta_1 -
 \alpha_{i} - \epsilon_i k(L_i)$.
 
 Choose $p_{i+1} \leq \bar{p}_2$ such that $p_{i+1}(\sigma)[L_{i+1}] =
 \bar{p}_2(\sigma)[L_{i+1}]$ for each $\sigma\in D(S_{i+1})$ --- this
 implies that Condition~\ref{co:7-1} holds --- and for each $\tau\in
 \prod_{\sigma\in D(S_{i+1})}p_{i+1}(\sigma)[L_{i+1}]$ there is some
 $C_{i+1,\tau}$ and $k_{i+1,\tau}$ such that
 $$p_{i+1}\angbr{\tau}\forces{\Poset^\kappa}{ \check{C}_{i+1,\tau}=
   C\cap\Rationals[k_{i+1,\tau}] \AND C\setminus \check{C}_{i+1,\tau}
   \in\Cantor^{\infty}_{\epsilon_{i+1}}} .$$
 For each $t\in
 p_i(\bar{\sigma})[L_i]$ and $x\in {\mathcal
   S}_{p_i(\bar{\sigma})}(t^*)$ and $\tau \in \prod_{\sigma\in
   D(S_{i})\setminus \{\bar{\sigma}\}}p_{i}(\sigma)[L_{i}]$ let
 $\rho(t,x,\tau)$ be the unique element of $\prod_{\sigma\in
   D(S_{i+1})}p_{i+1}(\sigma)[L_{i+1}]$ such that
 $\rho(t,x,\tau)\supseteq \tau(\sigma)$ for each $\sigma \in
 D(S_i)\setminus \{\bar{\sigma}\}$ and $\rho(t,x,\tau)\supseteq {t^*}$
 and $\rho(t,x,\tau)(|t^*|) = x$.  Let $\rho^*(t,x,\tau) $ be the
 unique element of $\prod_{\sigma\in D(S_{i})}p_{i}(\sigma)[L_{i}]$
 such that $\rho^*(t,x,\tau) = \tau$ if $S_i(\bar{\sigma}) = 0$ and
 $\rho^*(t,x,\tau) \supseteq \tau$ and $\rho^*(t,x,\tau)(\bar{\sigma})
 = t$ if $S_i(\bar{\sigma}) > 0$.

 Note that $C_{i+1,\rho(t,x,\tau)}\setminus C_{i,\rho^*(t,x,\tau)}\in
 \Cantor_{\epsilon_i}$ by Condition~\ref{co:99}.  Define
 $F_{t}:{\mathcal S}_{p_{i+1}(\bar{\sigma})}(t^*) \to
 \Cantor_{\epsilon_ik(L_i)^i}$ by $$F_{t}(x) = \bigcup_{\tau \in
   \prod_{\sigma\in D(S_{i})\setminus
     \{\bar{\sigma}\}}p_{i}(\sigma)[L_{i}]}
 C_{i+1,\rho(t,x,\tau)}\setminus C_{i,\rho^*(t,x,\tau)}
 $$
 and observe that $F_{t}(x)\in \Cantor_{\Gamma(m)/2}$
 by~\ref{eq:r1} since $$\left|\prod_{\sigma\in
     D(S_{i})}p_{i}(\sigma)[L_{i}]\right| \leq k(L_i)^{|D(S_i)|}\leq
 k(L_i)^i .$$
 
 By~\ref{eq:r4} it follows that $$\square_{\Theta(m-1)}(\{z \in
 \Cantor : \nu(\angbr{z,F_{t}}) < m\}) < \Gamma(m) $$
 for each $t\in
 p_i(\bar{\sigma})[L_i]$.  Let $b^t_{i+1}$ be such that
 $\lambda(b^t_{i+1}) < \Theta(m-1)$ and such that there is $a^t_{i+1}
 \in \Cantor_{\Gamma(m)}$ such that
\begin{equation}
  \label{eq:blt}
[a^t_{i+1}]\supseteq \{z \in \Cantor :
\nu(\angbr{z,F_{t}}) < m\}\setminus [b_{i+1}^t] .  
\end{equation}
It follows from~\ref{eq:r3} and Induction Hypothesis~\ref{co:101} that
if $b_{i+1}$ is defined to be $b_i \cup \bigcup_{t\in
  p_i(\bar{\sigma})[L_i]}b_{i+1}^t$ then Condition~\ref{co:101} is
satisfied by $b_{i+1}$.  Similarly, if $a_{i+1}$ is defined to be $a_i
\cup \bigcup_{t\in p_i(\bar{\sigma})[L_i]}a_{i+1}^t$ then $a_{i+1} \in
\Cantor_{\alpha_{i+1}}$ where $\alpha_{i+1} = \alpha_i +
\Gamma(m)k(L_i) < \delta_1$ by~\ref{eq:r2} and Induction
Hypothesis~\ref{co:10}.

In order to verify that Condition~\ref{co:1300} holds let $z\in
\Cantor \setminus [a_{i+1}\cup b_{i+1}]$. Let $T_{z,\sigma,i+1} =
T_{z,\sigma,i}$ if $\sigma\neq \bar{\sigma}$ and let
$T_{z,\bar{\sigma},i+1}$ be the set of all $s\in
p_{i+1}(\bar{\sigma})[\leq L_{i+1}]$ such that if $t^*\subseteq s$
then $s(|t^*|)\in \angbr{z,F_{t}}$. In order to show that
Condition~\ref{co:11b} holds it suffices to show that
$$\nu({\mathcal S}_{T_{\bar{\sigma},z,i+1}}(t^*)) \geq m >
S_{i+1}(\bar{\sigma})$$
for any $t \in p_i(\bar{\sigma})[L_i]$.  For
any given $t$ this follows from the fact that $z\notin [ a^t_{i+1}\cup
b^t_{i+1}]$ and \ref{eq:blt}.

In order to show that Condition~\ref{co:11c} holds let $\tau \in
\prod_{\sigma \in D(S_{i+1})}T_{{\sigma},z,i}[ L_{i+1}]$.  Let
$\bar{\tau} = \tau\restriction (D(S_i)\setminus \{\bar{\sigma}\})$.
Let $t\in p_i(\bar{\sigma})[L_i]$ be such that $t\subseteq
\tau(\bar{\sigma})$ and note that the definition of
$T_{\bar{\sigma},z,i}$ guarantees that $\tau(\bar{\sigma})(|t^*|)\in
\angbr{z,F_t}$. In other words, $z\notin
[F_t(\tau(\bar{\sigma})(|t^*|))]\supseteq
[C_{i+1,\rho(t,\tau(\bar{\sigma})(|t^*|),\bar{\tau})}\setminus
C_{i,\rho^*(t,\tau(\bar{\sigma})(|t^*|),\bar{\tau})}]$.  Moreover,
since $\rho^*(t,\tau(\bar{\sigma})(|t^*|),\bar{\tau})\in
\prod_{\sigma\in D(S_i)}p_i(\sigma)[L_i]$ it follows that $z\notin
[C_{i,\rho^*(t,\tau(\bar{\sigma})(|t^*|),\bar{\tau})}]$ by the
Induction Hypothesis~\ref{co:11c}. Hence $z\notin [C_{i+1,\tau}]$.
\end{proof}

\begin{corol}\label{l:2}
  If $G\subseteq \Poset$ is generic over $V$ then $V[G]\models
  \lambda^r_\infty(V\cap \Cantor) = 1$.
\end{corol}
\begin{proof}Suppose not and that 
  $$p\forces{\Poset^\kappa}{C\in \Cantor^{\infty}_\delta \AND V\cap
    \Cantor \subseteq [C]} $$
  and $\delta < 1$.  Using Lemma~\ref{l:0}
  find $A\subseteq \Cantor$ and $B\subset \Cantor$ such that
  $\lambda^r_\infty(A) + \lambda(B) < 1$ and $q\leq p$ such that
  $q\notforces{\Poset^\kappa}{z\in [C]}$ for each $z \in
  \Cantor\setminus (A \cup B)$.  Choose any $z \in \Cantor\setminus (A
  \cup B)$ and $q'\leq q$ such that $q'\forces{\Poset^\kappa}{z\notin
    [C]}$.
\end{proof}
\section{Remarks and open questions}
One might expect that the methods developed here could be used to
prove that for any two reals $s$ and $t$ such that $0 < t < s< 1$ it
is consistent that sets of size $\aleph_1$ are null with respect to
$s$-dimensional Hausdorff measure but that this is not so for
$t$-dimensional Hausdorff measure. While this is true of most of the
argument there are some slippery spots. For example. the use of $Y_F$
in the proof of Lemma~\ref{finitelymanyfunctions} assumes the
$\sigma$-finiteness of Lebesgue measure. Lemma~\ref{manylp} might also
pose some challenges to generalization.  Let ${\mathfrak n}_s$ denote
the least cardinal of a set which is not null with respect to
$s$-dimensional Hausdorff measure Hence the following question remains
open:
\begin{quest}
  Is it consistent that $0 < t < s< 1$ and ${\mathfrak n}_t
  <{\mathfrak n}_s$?
\end{quest}
\begin{quest}
  Is it consistent that $0 < u < v < w <1$ and ${\mathfrak n}_u
  <{\mathfrak n}_v < {\mathfrak n}_w$?
\end{quest}
\begin{quest}
  How big can the cardinality of $\{{\mathfrak n}_s\}_{s\in (0,1)}$
  be?
\end{quest}
However, the main open problem in this area still remains
Question~\ref{Komjath}. It would be interesting to know what the
answer to this question is in the model described in \S\ref{s:tfpo}.

\bibliography{../myabbrev,../standard,../steprans}



\begin{thebibliography}{1}

\bibitem{fpl}
David Fremlin.
\newblock Problem sheet.
\newblock {\tt http://www.essex.ac.uk/maths/staff/fremlin/problems.htm}, 2002.

\bibitem{GoShMan93}
Martin Goldstern and Saharon Shelah.
\newblock Many simple cardinal invariants.
\newblock {\em Archive for Mathematical Logic}, 32:203--221, 1993.

\bibitem{MR96h:28006}
Pertti Mattila.
\newblock {\em Geometry of sets and measures in {E}uclidean spaces}, volume~44
  of {\em Cambridge Studies in Advanced Mathematics}.
\newblock Cambridge University Press, Cambridge, 1995.

\bibitem{MR1613600}
Andrzej Ros{\l}anowski and Saharon Shelah.
\newblock Norms on possibilities. {I}. {F}orcing with trees and creatures.
\newblock {\em Mem. Amer. Math. Soc.}

\bibitem{step.36}
T.~Salisbury and J.~Stepr\={a}ns.
\newblock {H}ausdorff measure and {L}ebesgue measure.
\newblock {\em Real Anal. Exchange}, 22:265--278, 1996/97.

\bibitem{step.37}
J.~Stepr\={a}ns.
\newblock Cardinal invariants associated with {H}ausdorff capacities.
\newblock In {\em Proceeding of the {B}{E}{S}{T} Conferences 1-3}, volume 192
  of {\em Contemporary Mathematics}, pages 157--184. AMS, 1995.

\bibitem{step.42}
J.~Stepr\={a}ns.
\newblock Unions of rectifiable curves and the dimension of banach spaces.
\newblock {\em J. Symbolic Logic}, 64(2):701--726, June 1999.

\end{thebibliography}
\end{document}